\documentclass[12pt]{article}
\usepackage{amssymb}
\newcommand{\wt}{\widetilde}
\newcommand{\R}{\mathbb R}
\newcommand{\C}{\mathbb C}

\begin{document}
\thispagestyle{empty}
\hfill
\vskip 3.8truecm
\centerline{{\large The Beckman-Quarles theorem for continuous}}
\centerline{{\large mappings from ${\R}^n$ to ${\C}^n$}}
\vskip 0.8truecm
\centerline{{\large Apoloniusz Tyszka}}
\vskip 3.0truecm
\begin{abstract}
Let $\varphi_n:{\C}^n \times {\C}^n \to \C$,
$\varphi_n((x_1,...,x_n),(y_1,...,y_n))=(x_1-y_1)^2+...+(x_n-y_n)^2$.
We say that $f:{\R}^n \to {\C}^n$ preserves distance $d~\geq~0$ if for
each $x,y \in {\R}^n$ $\varphi_n(x,y)=d^2$ implies $\varphi_n(f(x),f(y))=d^2$.
We prove that if $x,y \in {\R}^n$ ($n \geq 3$) and
$|x-y|=(\sqrt{2+2/n})^k \cdot (2/n)^l$ ($k,l$ are non-negative
integers) then there exists a finite set
$\{x,y\} \subseteq S_{xy} \subseteq {\R}^n$ such that each
unit-distance preserving mapping from
$S_{xy}$ to ${\C}^n$ preserves the distance
between $x$ and $y$. It implies that each continuous map from
${\R}^n$ to ${\C}^n$ ($n \geq 3$) preserving unit distance
preserves all distances.
\end{abstract}
\vskip 0.2truecm
\normalsize
\par
The classical Beckman-Quarles theorem states that each unit-distance
preserving mapping from ${\R}^n$ to ${\R}^n$ ($n \geq 2$) is an isometry,
see \cite{Beckman}, \cite{Benz} and~\cite{Lester}.
Author's discrete form of this theorem
(\cite{Aequationes},\cite{Yerevan}) states that if
$x,y \in {\R}^n$ ($n \geq 2$) and $|x-y|$ is an algebraic number then
there exists a finite set $\{x,y\} \subseteq S_{xy} \subseteq {\R}^n$
such that each unit-distance preserving mapping from $S_{xy}$ to
${\R}^n$ preserves the distance between $x$ and $y$.
\vskip 0.2truecm
\par
Let $\varphi_n:{\C}^n \times {\C}^n \to \C$,
$\varphi_n((x_1,...,x_n),(y_1,...,y_n))=(x_1-y_1)^2+...+(x_n-y_n)^2$.
We say that $f:{\R}^n \to {\C}^n$ preserves distance $d \geq 0$ if for
each $x,y \in {\R}^n$ $\varphi_n(x,y)=d^2$ implies $\varphi_n(f(x),f(y))=d^2$.
Our goal is to prove that each continuous map from ${\R}^n$ to
${\C}^n$ ($n \geq 3$) preserving unit distance preserves all distances.
It requires some technical propositions.
\vskip 0.2truecm
\par
{\bf Proposition~1} (cf. \cite{Blumenthal}, \cite{Borsuk}).
The points $c_{1}=(z_{1,1},...,z_{1,n}),...,
c_{n+1}=(z_{n+1,1},...,z_{n+1,n}) \in {\C}^n$ are affinely dependent
if and only if their Cayley-Menger determinant
$$
\det \left[
\begin{array}{ccccc}
 0     &  1                        &  1                       & ... & 1                          \\
 1     & \varphi_n(c_{1},c_{1})    & \varphi_n(c_{1},c_{2})   & ... & \varphi_n(c_{1},c_{n+1})   \\
 1     & \varphi_n(c_{2},c_{1})    & \varphi_n(c_{2},c_{2})   & ... & \varphi_n(c_{2},c_{n+1})   \\
...    & ...                       & ...  	              & ... & ...                        \\
 1     & \varphi_n(c_{n+1},c_{1})  & \varphi_n(c_{n+1},c_{2}) & ... & \varphi_n(c_{n+1},c_{n+1}) \\
\end{array}\;\right]
$$
\par
\noindent
equals $0$.
\par
\vskip 0.2truecm
{\it Proof.} It follows from the equality
$$
\left(
\det \left[
\begin{array}{ccccc}
 z_{1,1}   & z_{1,2}   & ... &  z_{1,n}  & 1   \\
 z_{2,1}   & z_{2,2}   & ... &  z_{2,n}  & 1   \\
 ...       & ...       & ... &  ...      & ... \\
 z_{n+1,1} & z_{n+1,2} & ... & z_{n+1,n} & 1   \\
\end{array}
\right] \right)^2=
$$
$$
\frac{(-1)^{n+1}}{2^{n}} \cdot
\det \left[
\begin{array}{ccccc}
 0     &  1                       & 1                        & ... &  1                         \\
 1     & \varphi_n(c_{1},c_{1})   & \varphi_n(c_{1},c_{2})   & ... & \varphi_n(c_{1},c_{n+1})   \\
 1     & \varphi_n(c_{2},c_{1})   & \varphi_n(c_{2},c_{2})   & ... & \varphi_n(c_{2},c_{n+1})   \\
...    & ...                      & ...                      & ... & ...                        \\
 1     & \varphi_n(c_{n+1},c_{1}) & \varphi_n(c_{n+1},c_{2}) & ... & \varphi_n(c_{n+1},c_{n+1}) \\
\end{array}
\right].
$$
\vskip 0.2truecm
\par
{\bf Proposition~2} (cf. \cite{Blumenthal}, \cite{Borsuk}).
For each points $c_{1},...,c_{n+2} \in {\C}^n$ their
Cayley-Menger determinant equals $0$ i.e.
$$
\det \left[
\begin{array}{ccccc}
 0     &  1                       &  1                       & ... & 1                          \\
 1     & \varphi_n(c_{1},c_{1})   & \varphi_n(c_{1},c_{2})   & ... & \varphi_n(c_{1},c_{n+2})   \\
 1     & \varphi_n(c_{2},c_{1})   & \varphi_n(c_{2},c_{2})   & ... & \varphi_n(c_{2},c_{n+2})   \\
...    & ...                      & ...  	             & ... & ...                        \\
 1     & \varphi_n(c_{n+2},c_{1}) & \varphi_n(c_{n+2},c_{2}) & ... & \varphi_n(c_{n+2},c_{n+2}) \\
\end{array}\;\right]
={\rm 0}.
$$
\vskip 0.3truecm
{\it Proof.} Assume that
$c_{1}=(z_{1,1},...,z_{1,n}),...,c_{n+2}=(z_{n+2,1},...,z_{n+2,n})$.
The points
$\wt{c}_{1}=(z_{1,1},...,z_{1,n},0),...,
\wt{c}_{n+2}=(z_{n+2,1},...,z_{n+2,n},0) \in {\C}^{n+1}$ are affinely
dependent. Since $\varphi_n(c_i,c_j)=\varphi_{n+1}(\wt{c}_{i},\wt{c}_{j})$
$(1 \leq i \leq j \leq n+2)$ the Cayley-Menger determinant
of points $c_{1},...,c_{n+2}$ is equal to the Cayley-Menger determinant
of points $\wt{c}_{1},...,\wt{c}_{n+2}$ which equals $0$ according to
Proposition 1.
\par
\vskip 0.3truecm
\par
From Proposition~1 we obtain the following Propositions~3a and 3b.
\vskip 0.3truecm
\par
{\bf Proposition~3a.} If $d>0$ and points $c_1,...,c_{n+1} \in {\C}^n$
satisfy $\varphi_n(c_i,c_j)=d^2$ ($1 \leq i<j \leq n+1$),
then points $c_{1},...,c_{n+1}$ are affinely independent.
\vskip 0.3truecm
\par
{\bf Proposition~3b.} If $d>0$ and points
$c_1,...,c_n, c \in {\C}^n$ satisfy
$\varphi_n(c_i,c_j)=(2+2/n) \cdot d^2$ ($1 \leq i < j \leq n$) and
$\varphi_n(c_i,c)=d^2$ ($1 \leq i \leq n$), then points
$c_{1},...,c_{n},c$ are affinely independent.
\vskip 0.3truecm
\par
{\bf Proposition~4} (cf. \cite{Borsuk} p. 127 in the real case).
If points $c_{0},c_{1},...,c_{n} \in {\C}^n$ are affinely
independent, $x,y \in {\C}^n$ and
$\varphi_n(x,c_{0})=\varphi_n(y,c_{0})$, 
$\varphi_n(x,c_{1})=\varphi_n(y,c_{1})$,...,$\varphi_n(x,c_{n})=\varphi_n(y,c_{n})$,
then $x=y$.
\vskip 0.2truecm
\par
{\it Proof.} Computing we obtain that the vector
$\overrightarrow{xy}:=[s_1,...,s_n]$ is perpendicular
to each of the $n$ linearly independent vectors
$\overrightarrow{c_{0}c_{i}}$ ($i=1,...,n$).
Thus the vector $\overrightarrow{xy}$ is perpendicular
to every linear combination of vectors 
$\overrightarrow{c_{0}c_{i}}$ ($i=1,...,n$). 
In particular, the vector
$\overrightarrow{xy}=[s_1,...,s_n]$ is perpendicular
to the vector $[\bar{s_1},...,\bar{s_n}]$, where $\bar{s_1},...,\bar{s_n}$
denote numbers conjugate to the numbers $s_1,...,s_n$. Therefore
$\overrightarrow{xy}=0$ and the proof is complete.
\vskip 0.2truecm
\par
{\bf Proposition~5.} For each $n \in \{3,4,5,...\}$ there exists
$k(n) \in \{0,1,2,...\}$ such that
$1/2 \leq (2/n) \cdot (\sqrt{2+2/n})^{k(n)}<1$.
\vskip 0.3truecm
\par
{\bf Proposition~6.} If $n \in \{3,4,5,...\}$ then the set
$\{(\sqrt{2+2/n})^k \cdot (2/n)^l: k,l \in \{0,1,2,...\}\}$ is a
dense subset of $(0,\infty)$.
\vskip 0.2truecm
\par
{\it Proof.} Equivalently, our statement says that the set
$\{{\rm ln}_{\sqrt{2+2/n}}((\sqrt{2+2/n})^k \cdot (2/n)^l):
k,l \in \{0,1,2,...\}\}$ is a dense
subset of $\R$. Computing we obtain the set
$\{k+l \cdot {\rm ln}_{\sqrt{2+2/n}}(2/n): k,l \in \{0,1,2,...\}\}$
which is a dense subset of $\R$ due to
Kronecker's theorem (\cite{Hardy}) because
${\rm ln}_{\sqrt{2+2/n}}{(2/n)}$
is irrational and negative.
\vskip 0.3truecm
\par
{\bf Theorem~1.} If $x,y \in {\R}^n$ ($n\geq 3$) and
$|x-y|=(\sqrt{2+2/n})^k \cdot (2/n)^l$ ($k,l$ are non-negative
integers), then
there exists a finite set $\{x,y\} \subseteq S_{xy} \subseteq {\R}^n$
such that each unit-distance preserving mapping from $S_{xy}$ to ${\C}^n$
preserves the distance between $x$ and $y$.
\vskip 0.3truecm
\par
{\it Proof.}
Let $D_{n}$ denote the set of all positive numbers $d$
with the following property:
\vskip 0.1truecm
\noindent
if $x,y \in {\R}^n$ and $|x-y|=d$ then there exists a finite set
$\{x,y\} \subseteq S_{xy} \subseteq {\R}^n$ such that any map
$f:S_{xy}\rightarrow {\C}^n$ that preserves unit distance
preserves also the distance between $x$ and $y$.
\vskip 0.2truecm
\par
Obviously $1 \in D_n$.
\vskip 0.2truecm
\par
{\bf Lemma~1.} If $n \in \{2,3,4,...\}$ and $d \in D_n$,
then $\sqrt{2+2/n} \cdot d \in D_n$.
\vskip 0.2truecm
\par
{\it Proof.}
Let $x,y \in {\R}^n$, $|x-y|=\sqrt{2+2/n} \cdot d$.
There exist points $p_{1},...,p_{n},
\wt{y},\wt{p}_{1},...,\wt{p}_{n} \in {\R}^n$ such that:
\\
$|x-p_{i}|=|y-p_{i}|=d$ ($1 \leq i \leq n$),
\\
$|p_{i}-p_{j}|=d$ ($1 \leq i<j \leq n$),
\\
$|x-\wt{y}|=\sqrt{2+2/n} \cdot d$,
\\
$|y-\wt{y}|=d$,
\\
$|x-\wt{p}_{i}|=|\wt{y}-\wt{p}_{i}|=d$ ($1 \leq i \leq n$),
\\
$|\wt{p}_{i}-\wt{p}_{j}|=d$ ($1 \leq i<j \leq n$).
\vskip 0.6truecm
\par
\noindent
Let
$$
S_{xy}:=
\bigcup_{i=1}^{n} S_{xp_{i}}
\cup
\bigcup_{i=1}^{n} S_{yp_{i}}
\cup
\bigcup_{1\le i<j \leq n} S_{p_{i}p_{j}}
\cup
\bigcup_{i=1}^{n} S_{x\wt{p}_{i}}
\cup
\bigcup_{i=1}^{n} S_{\wt{y}\wt{p}_{i}}
\cup
\bigcup_{1\le i<j \leq n} S_{\wt{p}_{i}\wt{p}_{j}}
\cup
S_{y\wt{y}}
$$
and $f: S_{xy} \rightarrow {\C}^n$
preserves unit distance. Since
$$
S_{xy} \supseteq
\bigcup_{i=1}^{n}S_{xp_{i}}
\cup
\bigcup_ {i=1}^{n}S_{yp_{i}}
\cup
\bigcup_{1 \leq i<j \leq n} S_{p_{i}p_{j}}
$$
\noindent
we conclude that $f$ preserves the distances between
$x$ and $p_{i}$ ($1\leq i \leq n$), $y$ and $p_{i}$ ($1\leq i\leq n$),
and all distances between $p_{i}$ and $p_{j}$ ($1\leq i<j\leq n$).
Hence for all $1 \leq i \leq n$ $\varphi_n(f(x),f(p_i))=\varphi_n(f(y),f(p_i))=d^2$
and for all $1 \leq i<j \leq n$ $\varphi_n(f(p_i),f(p_j))=d^2$.
Since $S_{xy} \supseteq S_{y\wt{y}}$ we conclude that
$\varphi_n(f(y),f(\wt{y}))=d^2$.
By Proposition~2 the Cayley-Menger determinant of points
$x$, $p_1,...,p_n, y$ equals $0$ i.e.
$$
\det \left[
\begin{array}{ccccccc}
 0     &  1                     &  1                       & ... &  1                      &  1                       \\
 1     & \varphi_n(f(x),f(x))   & \varphi_n(f(x),f(p_1))   & ... & \varphi_n(f(x),f(p_n))    & \varphi_n(f(x),f(y))   \\
 1     & \varphi_n(f(p_1),f(x)) & \varphi_n(f(p_1),f(p_1)) & ... & \varphi_n(f(p_1),f(p_n))  & \varphi_n(f(p_1),f(y)) \\
...    & ...                    & ...                      & ... & ...                     & ...                      \\
 1     & \varphi_n(f(p_n),f(x)) & \varphi_n(f(p_n),f(p_1)) & ... & \varphi_n(f(p_n),f(p_n))  & \varphi_n(f(p_n),f(y)) \\
 1     & \varphi_n(f(y),f(x))   & \varphi_n(f(y),f(p_1))   & ... & \varphi_n(f(y),f(p_n))    & \varphi_n(f(y),f(y))   \\
\end{array}
\right]
=0.
$$
Denoting $t=\varphi_n(f(x),f(y))$ we obtain
$$
\det \left[
\begin{array}{ccccccc}
 0     &  1     &  1     & ... &  1     &  1     \\
 1     &  0     & d^2    & ... & d^2    &  t     \\
 1     & d^2    &  0     & ... & d^2    & d^2    \\
...    & ...    & ...    & ... & ...    & ...    \\
 1     & d^2    & d^2    & ... &  0     & d^2    \\
 1     & t      & d^2    & ... & d^2    &  0     \\
\end{array}
\right]
=0.
$$
\par
\noindent
Computing this determinant we obtain
$$
(-1)^{n-1} \cdot d^{2n-2} \cdot t \cdot
\left( n \cdot t - (2n+2) \cdot d^2 \right)=0.
$$
Therefore
$$t=\varphi_n(f(x),f(y))=\varphi_n(f(y),f(x))=(\sqrt{2+2/n} \cdot d)^2$$
or $$t=\varphi_n(f(x),f(y))=\varphi_n(f(y),f(x))=0.$$
Analogously we may prove that
$$\varphi_n(f(x),f(\wt{y}))=\varphi_n(f(\wt{y}),f(x))=
(\sqrt{2+2/n} \cdot d)^2$$
or $$\varphi_n(f(x),f(\wt{y}))=\varphi_n(f(\wt{y}),f(x))=0.$$
If $t=0$ then the points $f(x)$ and $f(y)$ satisfy:
\\
\\
\centerline{$\varphi_n(f(x),f(x))=0=\varphi_n(f(y),f(x))$,}
\\
\centerline{$\varphi_n(f(x),f(p_1))=d^2=\varphi_n(f(y),f(p_1))$,}
\\
\centerline{$...$}
\\
\centerline{$\varphi_n(f(x),f(p_n))=d^2=\varphi_n(f(y),f(p_n))$.}
\\
\par
\noindent
By Proposition~3a the points
$f(x),f(p_1),...,f(p_n)$ are affinely independent.
Therefore by Proposition~4 $f(x)=f(y)$ and consequently
$$
d^2=\varphi_n(f(y),f(\wt{y}))=\varphi_n(f(x),f(\wt{y}))
\in \{(\sqrt{2+2/n} \cdot d)^2,~0 \}.
$$
Since $d^2 \neq (\sqrt{2+2/n} \cdot d)^2$ and $d^2 \neq 0$ we conclude that
the case $t=0$ cannot occur. This completes the proof of Lemma~1.
\vskip 0.4 truecm
\par
By Lemma~1 if $d \in D_n$ ($n \geq 2$) then all distances
$(\sqrt{2+2/n})^m \cdot d$ ($m=0,1,2,...$)
belong to $D_n$.
\vskip 0.4truecm
\par
{\bf Lemma~2.} If $n \in \{3,4,5,...\}$, $\varepsilon, d \in D_n$
and $\varepsilon/2 \leq (2/n) \cdot d \neq \varepsilon$,
then $(2/n) \cdot d \in D_n$.
\vskip 0.2truecm
\par
{\it Proof.} Let $x,y \in {\R}^n$, $|x-y|=(2/n) \cdot d$.
There exist points $p_{1},...,p_{n},\wt{y}$,
$\wt{p}_{1},...,\wt{p}_{n} \in {\R}^n$
such that:
\par
\noindent
$|x-p_{i}|=|y-p_{i}|=d$ ($1 \leq i \leq n$),
\\
$|p_{i}-p_{j}| = \sqrt{2+2/n} \cdot d$ ($1 \leq i \leq n$),
\\
$|x-\wt{y}|=(2/n) \cdot d$,
\\
$|y-\wt{y}|=\varepsilon$,
\\
$|x-\wt{p}_{i}|=|\wt{y}-\wt{p}_{i}|=d$ ($1 \leq i \leq n$),
\\
$|\wt{p}_{i}-\wt{p}_{j}|=\sqrt{2+2/n} \cdot d$ ($1 \leq i <j \leq n$).
\\
\\
Let
$$
S_{xy}:=
\bigcup_{i=1}^{n} S_{xp_{i}}
\cup
\bigcup_{i=1}^{n} S_{yp_{i}}
\cup
\bigcup_{1\le i<j \leq n} S_{p_{i}p_{j}}
\cup
\bigcup_{i=1}^{n} S_{x\wt{p}_{i}}
\cup
\bigcup_{i=1}^{n} S_{\wt{y}\wt{p}_{i}}
\cup
\bigcup_{1\le i<j \leq n} S_{\wt{p}_{i}\wt{p}_{j}}
\cup
S_{y\wt{y}}
$$
and $f: S_{xy} \rightarrow {\C}^n$ preserves unit distance. Since
$$
S_{xy} \supseteq
\bigcup_{i=1}^{n}S_{xp_{i}}
\cup
\bigcup_ {i=1}^{n}S_{yp_{i}}
\cup
\bigcup_{1 \leq i<j \leq n} S_{p_{i}p_{j}}
$$
\noindent
we conclude that $f$ preserves the distances between
$x$ and $p_{i}$ ($1\leq i \leq n$), $y$ and $p_{i}$ ($1\leq i\leq n$),
and all distances between $p_{i}$ and $p_{j}$ ($1\leq i<j\leq n$).
Hence for all $1 \leq i \leq n$ $\varphi_n(f(x),f(p_i))=\varphi_n(f(y),f(p_i))=d^2$
and for all $1 \leq i<j \leq n$ $\varphi_n(f(p_i),f(p_j))=(\sqrt{2+2/n} \cdot d)^2$.
Since $S_{xy} \supseteq S_{y\wt{y}}$ we conclude that
$\varphi_n(f(y),f(\wt{y}))=\varepsilon^2$.
By Proposition 2~the Cayley-Menger determinant of points
$x$, $p_1,...,p_n, y$ equals $0$ i.e.
$$
\det \left[
\begin{array}{ccccccc}
 0     &  1                     &  1                       & ... &  1                        &  1                     \\
 1     & \varphi_n(f(x),f(x))   & \varphi_n(f(x),f(p_1))   & ... & \varphi_n(f(x),f(p_n))    & \varphi_n(f(x),f(y))   \\
 1     & \varphi_n(f(p_1),f(x)) & \varphi_n(f(p_1),f(p_1)) & ... & \varphi_n(f(p_1),f(p_n))  & \varphi_n(f(p_1),f(y)) \\
...    & ...                    &  ...                     & ... &  ...                      &  ...                   \\
 1     & \varphi_n(f(p_n),f(x)) & \varphi_n(f(p_n),f(p_1)) & ... & \varphi_n(f(p_n),f(p_n))  & \varphi_n(f(p_n),f(y)) \\
 1     & \varphi_n(f(y),f(x))   & \varphi_n(f(y),f(p_1))   & ... & \varphi_n(f(y),f(p_n))    & \varphi_n(f(y),f(y))   \\
\end{array}
\right]
{\rm =0.}
$$
Denoting $t=\varphi_n(f(x),f(y))$ we obtain
$$
\det \left[
\begin{array}{ccccccc}
 0     &  1     &  1                & ... &  1                &  1      \\
 1     &  0     & d^2               & ... & d^2               &  t      \\
 1     & d^2    &  0                & ... & (2+2/n) \cdot d^2 & d^2     \\
...    & ...    & ...               & ... & ...               &  ...    \\
 1     & d^2    & (2+2/n) \cdot d^2 & ... &  0                & d^2     \\
 1     &  t     & d^2               & ... & d^2               &  0      \\
\end{array}
\right]
=0.
$$
\par
\noindent
Computing this determinant we obtain 
$$
\frac{(-2n-2)^{n-1}}{n^n} \cdot d^{2n-2} \cdot t \cdot
\left( n^2 \cdot t - 4d^2 \right)
=0.
$$
Therefore $$t=\varphi_n(f(x),f(y))=\varphi_n(f(y),f(x))=((2/n) \cdot d)^2$$
or $$t=\varphi_n(f(x),f(y))=\varphi_n(f(y),f(x))=0.$$ Analogously we may prove
that $$\varphi_n(f(x),f(\wt{y}))=\varphi_n(f(\wt{y}),f(x))=
((2/n)\cdot d)^2$$ or $$\varphi_n(f(x),f(\wt{y}))=\varphi_n(f(\wt{y}),f(x))=0.$$
If $t=0$ then the points $f(x)$ and $f(y)$ satisfy:
\\
\\
\centerline{$\varphi_n(f(x),f(x))=0=\varphi_n(f(y),f(x))$,}
\\
\centerline{$\varphi_n(f(x),f(p_1)=d^2=\varphi_n(f(y),f(p_1)$,}
\\
\centerline{$...$}
\\
\centerline{$\varphi_n(f(x),f(p_n))=d^2=\varphi_n(f(y),f(p_n))$.}
\\
\par
\noindent
By Proposition~3b the points
$f(x),f(p_1),...,f(p_n)$ are affinely independent.
Therefore by Proposition~4 $f(x)=f(y)$ and consequently
$$
\varepsilon^2=
\varphi_n(f(y),f(\wt{y}))=\varphi_n(f(x),f(\wt{y}))
\in \{((2/n) \cdot d)^2,~0\}.
$$
Since $\varepsilon^2 \neq ((2/n) \cdot d)^2$ and $\varepsilon^2 \neq 0$
we conclude that the case $t=0$ cannot occur.
This completes the proof of Lemma~2.
\vskip 0.4truecm
\par
By Proposition~5 for each $n \in \{3,4,5,...\}$
there exists $k(n) \in \{0,1,2,...\}$ such that if
$r \in D_n$ then
$\varepsilon:=r$ and $d:=(\sqrt{2+2/n})^{k(n)} \cdot r$
satisfy assumptions of Lemma~2 and
moreover $\rho(n):=(2/n) \cdot (\sqrt{2+2/n})^{k(n)} < 1$.
Therefore, by Lemma~2 if $r \in D_n$ then $\rho(n) \cdot r \in D_n$.
Since $1 \in D_{n}$ we conclude that all distances $\rho(n)^m$ ($m=0,1,2,...$)
belong to $D_n$. For each $d \in D_n$ there exists
$m \in \{0,1,2,...\}$ such that $\varepsilon:=\rho(n)^m < (2/n) \cdot d$.
Applying Lemma~2 for such $\varepsilon$
we obtain the following lemma.
\vskip 0.4truecm
\par
{\bf Lemma~3.} If $n \in \{3,4,5,...\}$ and $d \in D_n$, then
$(2/n) \cdot d \in D_n$.
\vskip 0.2truecm
\par
From Lemma~1 and Lemma~3 we obtain that for each non-negative
integers $k,l~~(\sqrt{2+2/n})^k \cdot (2/n)^l \in D_n$.
This completes the proof of Theorem~1.
\vskip 0.3truecm
\par
As a corollary of Theorem~1 and Proposition~6 we obtain our main theorem.
\vskip 0.2truecm
\par
{\bf Theorem~2.} Each continuous map from ${\R}^n$ to ${\C}^n$ ($n \geq 3$)
preserving unit distance preserves all distances.
\vskip 0.2truecm
\par
{\bf Remark.} By an endomorphism of $\C$ we understand any map
$f:\C \to \C$ satisfying:
\par
\centerline
{$\forall x,y \in \C ~~f(x+y)=f(x)+f(y)$,}
\centerline{$\forall x,y \in \C ~~f(x \cdot y)=f(x) \cdot f(y)$,}
\centerline{$f(0)=0$,}
\centerline{$f(1)=1$.}
\par
\noindent
If $f:\C \to \C$ is an endomorphism then
$(f_{|\R},...,f_{|\R}): {\R}^n \to {\C}^n$ preserves unit distance.
Bijective endomorphisms are called automorphisms.
There are two trivial automorphisms of $\C$:
identity and conjugation. It is known
that there exist non-trivial automorphisms of $\C$
and each such automorphism $f:\C \to \C$ satisfies:
$\exists^{x \in \R}_{x \neq 0} f(x) \not\in \R$ (\cite{Kuczma}). From this 
$\varphi_n((0,0,...,0),(x,0,...,0))=|x|^2$
and
$\varphi_n((f(0),f(0),...,f(0)),(f(x),f(0),...,f(0)))=(f(x))^2 \neq |x|^2$.
Therefore $(f_{|\R},...,f_{|\R})$ preserves unit distance, but does not
preserve the distance $|x|$~$>$~$0$.
\par

Technical Faculty\\
Hugo Ko{\l}{\l}\c{a}taj University\\
Balicka 104, 30-149 Krak\'ow, Poland\\
E-mail: {\it rttyszka@cyf-kr.edu.pl}\\
\end{document}